# Some Useful Collective Properties of Bessel, Marcum Q-Functions and Laguerre Polynomials

Hakan ÖZTÜRK[1*], Fikret ANLI[2] and Abdelouahab KADEM[3]

[1]Department of Physics, Faculty of Arts and Sciences, Osmaniye Korkut Ata University, 80000, Osmaniye, Turkey
[2]Department of Physics, Faculty of Arts and Sciences, Kahramanmaraş Sütçü Imam University, 46040, Kahramanmaraş, Turkey
[3]Department of Mathematics, Faculty of Sciences, Ferhat Abbas University, LMFN laboratory, Setif, Algeria

[*]Corresponding author (H. Öztürk)

**Abstract**

Special functions have been used widely in many problems of applied sciences. However, there are considerable numbers of problems in which exact solutions could not be achieved because of undefined sums or integrals involving special functions. These handicaps force researchers to seek new properties of special functions. Many problems that could not be solved so far would be solved by means of these efforts. Therefore in this article, we derived some useful properties and interrelations of each others of Bessel functions, Marcum Q-functions and Laguerre polynomials.

Funding: The authors declare that no funds, grants, or other support were received during the preparation of this manuscript.
Competing interest: The authors have no relevant financial or non-financial interests to disclose.

## 1. Introduction

Bessel functions are the solutions of a second order linear differential equation and they are particularly important for the problems of wave propagation and electrostatic potentials. These functions are not restricted to these areas; they are also frequently used in many of the problems such as inviscid rotational flows, heat conduction, diffusion problems, radial Schrödinger equation in cylindrical or spherical coordinates [1-3].

Similarly in quantum mechanics, the Schrödinger wave equation for the hydrogen-like atom can be solved exactly by the separation of variables in spherical coordinates. The radial parts of the wave functions presented in solution include the Laguerre polynomials. These polynomials are also used in describing the static Wigner functions of oscillatory systems in phase space. They are further mentioned in the quantum mechanics of the potential energy of a diatomic molecule and of the 3D isotropic harmonic oscillator [4].

The generalized Marcum Q-function is an important function used in radar detection and communications [5, 6] and also in statistics and probability theory, where they are called non-central chi-square or non-central gamma cumulative distributions [7]. Detailed information on Marcum Q-functions can be found in Ref. [8, 9].

A brief recallment about the basic properties of Bessel functions, Laguerre polynomials and Marcum Q-functions are given in section 2. Besides some advanced properties of these functions, a number of important and new relations between them are derived in section 3. Some comments are also given in conclusion section. During the derivation of the equations an explicit formalism is preferred rather than a closed representation for understanding the expressions better. This study is an extension of the previous studies whose are related with the new properties of the Bessel, Marcum Q-functions and Laguerre polynomials which are frequently used in the solutions of the problems vitally important in science and engineering.

## 2. Mathematical preliminaries

Before we start our derivation of the equations including new relations, it is better to give some existent properties of the Bessel functions, Laguerre polynomials and Marcum Q-functions, explicitly. The series expansion, integral representation and generating function of modified Bessel functions of the first kind for integer $n$ are given as, respectively [10];

$$I_n(x) = \sum_{k=0}^{\infty} \frac{(x/2)^{n+2k}}{k!(n+k)!} \tag{1}$$

$$I_n(x) = \frac{(2x)^n n!}{\pi (2n)!} \int_{-1}^{1} (1-\mu^2)^{n-1/2} \exp(-x\mu) d\mu, \quad n > -1/2 \tag{2}$$

$$\sum_{n=-\infty}^{\infty} t^n I_n(x) = \exp\left(\frac{x}{2}(t+1/t)\right), \quad I_{-n}(x) = I_n(x), \quad -1 \leq t \leq t \tag{3}$$

Laguerre polynomials are one of the most widely used among the classical orthogonal polynomials in physics and applied sciences including mathematics. Series expansion, integral representation and generating function of them are given as, respectively [11, 12];

$$L_n(x) = \sum_{k=0}^{n} \frac{(-1)^k n! x^k}{(k!)^2 (n-k)!} \tag{4}$$

$$L_n(x) = \frac{\exp(x)}{n!} \int_0^{\infty} u^n \exp(-u) J_0(2\sqrt{xu}) du \tag{5}$$

$$\sum_{n=0}^{\infty} t^n L_n(x) = \frac{1}{1-t} \exp\left(-\frac{xt}{1-t}\right), \quad -1 < t < 1 \tag{6}$$

where $J_0(\cdot)$ is the zeroth order Bessel functions of the first kind. Generalized Marcum Q-function is defined as in original notation [8],

$$Q_M(\alpha,\beta) = \int_\beta^\infty x\left(\tfrac{x}{\alpha}\right)^{M-1} \exp\left[-\tfrac{(x^2+\alpha^2)}{2}\right] I_{M-1}(\alpha x)\,dx \tag{7}$$

This integral representation of Marcum Q-functions satisfies the recurrence relations and formal Neumann series expansions, respectively [8];

$$Q_M(\alpha,\beta) - Q_{M-1}(\alpha,\beta) = \left(\tfrac{\beta}{\alpha}\right)^{M-1} \exp\left[-\tfrac{(\beta^2+\alpha^2)}{2}\right] I_{M-1}(\alpha\beta) \tag{8}$$

$$Q_M(\alpha,\beta) = 1 - \exp\left[-\tfrac{(\beta^2+\alpha^2)}{2}\right] \sum_{n=M}^\infty \left(\tfrac{\beta}{\alpha}\right)^n I_n(\alpha\beta) \tag{9}$$

where $I_M(\cdot)$ is the $M$th order modified Bessel functions of the first kind.

## 3. Properties and connections

In the first part of this study, we multiply Eq. (1) by $t^n$ and sum the resulting equation over $n$ from zero to infinity,

$$S(x,t) = \sum_{n=0}^\infty t^n I_n(x) = \exp\left[\tfrac{x}{2}\left(t+\tfrac{1}{t}\right)\right] - \exp\left(\tfrac{xt}{2}\right) \sum_{k=1}^\infty \tfrac{[x/(2t)]^k}{k!(k-1)!} \Gamma(k, xt/2) \tag{10}$$

where $\Gamma(k, xt/2)$ is the incomplete gamma function defined as,

$$\Gamma(k, xt/2) = \int_{xt/2}^\infty u^{k-1} \exp(-u)\,du. \tag{11}$$

Since the summation on the right hand side of Eq. (10) is unknown, it is unlikely to represent Eq. (10) in an explicit form. Therefore, we need to use other mathematical tools such as converting it to a differential equation given below;

$$\frac{\partial S(x,t)}{\partial x} - \frac{1}{2}\left(t+\tfrac{1}{t}\right) S(x,t) = -\frac{1}{2t} I_0(x) + \frac{1}{2} I_1(x) \tag{12}$$

where $I_0(x)$ and $I_1(x)$ are zeroth and first order modified Bessel functions of the first kind, respectively. For the solution of the partial differential equation given in Eq. (12), we get;

$$S(x,t) = \exp\left(\tfrac{x(t+1/t)}{2}\right)\left\{-\tfrac{1}{2t}\int I_0(x)\exp\left(\tfrac{-x(t+1/t)}{2}\right)dx + \tfrac{1}{2}\int I_1(x)\exp\left(\tfrac{-x(t+1/t)}{2}\right)dx + C_1(t)\right\} \tag{13}$$

Up to now, integrals in Eq. (13) were unknown, that is, analytic results of these integrals were defined in terms of known elementary and special functions. However we are not hopeless: By setting $M = 0$ and coordinating the variables as $\beta/\alpha = t$, $\alpha\beta = x$ in Eq. (9), we get,

$$S(x,t) = \sum_{n=0}^{\infty} t^n I_n(x) = \exp\left[\tfrac{x}{2}\left(t+\tfrac{1}{t}\right)\right]\left(1 - Q_0(\sqrt{x/t}, \sqrt{xt})\right) \tag{14}$$

Then, the integrals in Eq. (13) can now be evaluated as,

$$\int I_0(x)\exp\left(\tfrac{-x(t+1/t)}{2}\right)dx = \frac{4t}{1-t^2}Q_0(\sqrt{x/t},\sqrt{xt}) + \frac{2t}{1-t^2}I_0(x)\exp\left(\tfrac{-x(t+1/t)}{2}\right) + \frac{4t}{1-t^2}(C_1(t)-1) \tag{15}$$

$$\int I_1(x)\exp\left(\tfrac{-x(t+1/t)}{2}\right)dx = \frac{2(1+t^2)}{1-t^2}Q_0(\sqrt{x/t},\sqrt{xt}) + \frac{2}{1-t^2}I_0(x)\exp\left(\tfrac{-x(t+1/t)}{2}\right) + \frac{4}{1-t^2}(C_1(t)-1) \tag{16}$$

where $C_1(t)$ is the integration constant which can be determined by using Eq. (14). Setting $x = 0$ in Eq. (14), we have $S(0, t) = 1$ and $Q_0(0,0) = 0$. Then using these identities in Eq. (15) or Eq. (16) we find $C_1(t) = 1/2$. An alternative way to determine the integration constant $C_1(t)$ is to use the series expansion of Eq. (15) and Eq. (16) with respect to variable $x$. It is seen that Eq. (14) satisfies Eq. (12) exactly. Note that Eq. (14) is a key equation to derive the results of the integrals in Eq. (13). Eq. (14) has also another property: Generating function of modified Bessel functions of the first kind is given in Eq. (3) for more general definition, that is, summation is from $-\infty$ to $+\infty$ for integer $n$. However, Eq. (14) contains half space summation for $n$, then we refer to Eq. (14) as *the main generating function* of modified Bessel functions of the first kind.

In the second part of this study, a summation over Laguerre polynomials which does not exist in literature will be derived explicitly. Unfortunately, it is not known yet how to write Marcum Q-functions of zeroth order $Q_0(\alpha, \beta)$ in terms of known elementary functions. If it was known we could derive all the $Q_n(\alpha, \beta)$ functions using Eq. (8). By using Eq. (3) and Eq. (14) we can easily deduce that

$$Q_0(\sqrt{x/t},\sqrt{xt}) + Q_0(\sqrt{xt},\sqrt{x/t}) = 1 - \exp\left[-\tfrac{x}{2}\left(t+\tfrac{1}{t}\right)\right]I_0(x) \tag{17}$$

A different integral representation for $Q_0(\sqrt{x/t},\sqrt{xt})$ can be obtained by multiplying Eq. (2) by $t^n$ and by summing the resulting equation over $n$ from zero to infinity,

$$\begin{aligned}S(x,t) &= \sum_{n=0}^{\infty} t^n I_n(x) \\ &= I_0(x) + \sqrt{\tfrac{xt}{2\pi}}\exp\left[\tfrac{x}{2}\left(t+\tfrac{1}{t}\right)\right]\int_{-1}^{1}\tfrac{\mu}{\sqrt{1-\mu^2}}\exp\left[-\tfrac{xt}{2}(1-\mu^2)\right]erf\left(\sqrt{\tfrac{xt}{2}}\mu + \sqrt{\tfrac{x}{2t}}\right)d\mu\end{aligned} \tag{18}$$

Using Eq. (14) and Eq. (18) we get,

$$\begin{aligned}&\int_{-1}^{1}\tfrac{\mu}{\sqrt{1-\mu^2}}\exp\left[-\tfrac{xt}{2}(1-\mu^2)\right]erf\left(\sqrt{\tfrac{xt}{2}}\mu + \sqrt{\tfrac{x}{2t}}\right)d\mu \\ &= \sqrt{\tfrac{2\pi}{xt}}\left\{1 - Q_0(\sqrt{x/t},\sqrt{xt}) - I_0(x)\exp\left[-\tfrac{x}{2}\left(t+\tfrac{1}{t}\right)\right]\right\}\end{aligned} \tag{19}$$

where $erf(\cdot)$ is the error function defined as,

$$erf(z) = \tfrac{2}{\sqrt{\pi}} \int_0^z \exp(-u^2) \, du \qquad (20)$$

One of the most noticeable relations between modified Bessel functions of the first kind and Laguerre polynomials is given in Ref. [11-13] as,

$$\sum_{n=0}^{\infty} \frac{(-1)^n t^n L_n(x)}{n!} = \exp(-t) I_0(2\sqrt{xt}) \qquad (21)$$

It is possible to prove this relation by multiplying Eq. (5) by $(-1)^n t^n / n!$ and summing the resulting equation over $n$ from zero to infinity,

$$\sum_{n=0}^{\infty} \frac{(-1)^n t^n L_n(x)}{n!} = \exp(x) \int_0^{\infty} \exp(-u) J_0(2\sqrt{tu}) J_0(2\sqrt{xu}) \, du = \exp(x) \{ \exp(-x-t) I_0(2\sqrt{xt}) \} \qquad (22)$$

In order to get a more general expression, let us first multiply Eq. (5) by $(-1)^n t^n / (n+m)!$ and then sum the resulting equation over $n$ from zero to infinity; we get,

$$\sum_{n=0}^{\infty} \frac{(-1)^n t^n L_n(x)}{(n+m)!} = \exp(x) \int_0^{\infty} \frac{\exp(-u)}{(tu)^{m/2}} J_m(2\sqrt{tu}) J_0(2\sqrt{xu}) \, du \qquad (23)$$

If Eq. (23) is rearranged for $m = 1$, we reobtain the results of Pent [14] and Andras et al. [15],

$$\sum_{n=0}^{\infty} \frac{(-1)^n t^n L_n(x)}{(n+1)!} = \exp(x) \left\{ \frac{\exp(-x-t)}{t} \sum_{n=1}^{\infty} \left( \tfrac{t}{x} \right)^{n/2} I_n(2\sqrt{xt}) \right\} = \frac{\exp(-t)}{t} \left( 1 - Q_1(\sqrt{2x}, \sqrt{2t}) \right) \qquad (24)$$

Then, by evaluating the integral in Eq. (23), we get the following general formulation

$$\sum_{n=0}^{\infty} \frac{(-1)^n t^n L_n(x)}{(n+m)!} = \exp(x) \left\{ \frac{\exp(-x-t)}{t^m} \sum_{n=m}^{\infty} \frac{(n-1)!}{(m-1)!(n-m)!} \left( \tfrac{t}{x} \right)^{n/2} I_n(2\sqrt{xt}) \right\}, \quad m \geq 1 \qquad (25)$$

## 4. Conclusion

In this paper, all derivations are checked by means of one of the mathematical software (maple) and several unknown collective properties of Laguerre polynomials, Bessel and Marcum Q-functions are derived explicitly. Eq. (14), Eq. (15), Eq. (16) and Eq. (25) are new novelties in mathematical sciences. Unfortunately, we do not know how to write the Marcum Q-functions of zeroth order $Q_0(\alpha, \beta)$ in terms of the well-known elementary functions. However, we will give a few meaningful comments for it:

For $t = 1$, then we get from Eq. (14),

$$Q_0(\sqrt{x}, \sqrt{x}) = [1 - \exp(-x) I_0(x)] / 2 \qquad (26)$$

and from Eq. (19),

$$\int_{-1}^{1} \frac{\mu}{\sqrt{1-\mu^2}} \exp\left[-\tfrac{x}{2}(1-\mu^2)\right] \operatorname{erf}\left(\sqrt{\tfrac{x}{2}}(\mu+1)\right) d\mu = \sqrt{\tfrac{\pi}{2x}}\left(1-\exp(-x)I_0(x)\right) = \sqrt{\tfrac{2\pi}{x}} Q_0(\sqrt{x},\sqrt{x}) \quad (27)$$

For $t = -1$, then we get from Eq. (16),

$$Q_0(i\sqrt{x}, i\sqrt{x}) = \left[1 - \exp(x)I_0(x)\right]/2 \quad (28)$$

and from Eq. (19),

$$\int_{-1}^{1} \frac{\mu}{\sqrt{1-\mu^2}} \exp\left[\tfrac{x}{2}(1-\mu^2)\right] \operatorname{erf}\left(i\sqrt{\tfrac{x}{2}}(\mu+1)\right) d\mu = i\sqrt{\tfrac{\pi}{2x}}\left(1-\exp(x)I_0(x)\right) = i\sqrt{\tfrac{2\pi}{x}} Q_0(i\sqrt{x}, i\sqrt{x}) \quad (29)$$

where $i = \sqrt{-1}$.

Note that all the $Q_M(\alpha, \beta)$ functions are even functions of the variables $\alpha$ and $\beta$. Some limit values of $Q_0(\alpha, \beta)$ and $Q_M(\alpha, \beta)$, $M \geq 1$ functions can also be given as, respectively;

$$Q_0(\alpha, \beta) = \begin{cases} 0, & \alpha = 0 \\ 1, & \alpha = \pm\infty \\ 1 - \exp(-\alpha^2/2), & \beta = 0 \\ 1, & \beta = \pm\infty \end{cases} \quad (30)$$

$$Q_M(\alpha, \beta) = \begin{cases} \frac{\Gamma(M, \beta^2/2)}{(M-1)!}, & \alpha = 0 \\ 1, & \alpha = \pm\infty \\ 1, & \beta = 0 \\ 1, & \beta = \pm\infty \end{cases} \quad (31)$$

**References**

[1]. Abramowitz, M., Stegun, I.A.: Handbook of mathematical functions with formulas, graphs, and mathematical tables. 9th Ed., New York: Dover (1972).
[2]. Jimenez, N., Pico, R., Sanchez-Morcillo, V., Romero-Garcia V., Garcia-Raffi, L.M., Staliunas, K.: Formation of high-order acoustic Bessel beams by spiral diffraction gratings, Phys. Rev. E 94(5): 053004 (2016). https://doi.org/10.1103/PhysRevE.94.053004
[3]. Temme, N.M.: Special functions: An introduction to the classical functions of mathematical physics, 2nd print ed., New York, Wiley (1996).
[4]. Schatz, G.C.: Ratner MA. Quantum mechanics in chemistry, 0-13-895491-7, NJ: Prentice Hall (2001).
[5]. Marcum, J.I.: A statistical theory of target detection by pulsed radar, IRE Trans. Inform. Theory, 6(2), 59-267 (1960). DOI: 10.1109/TIT.1960.1057560

[6]. Rice, S.O.: Uniform asymptotic expansions for saddle point integrals-application to a probability distribution occurring in noise theory, Bell. System Tech. J., 47, 1971-2013 (1968). https://doi.org/10.1002/j.1538-7305.1968.tb01099.x

[7]. Ross, A.H.M.: Algorithm for calculating the noncentral chi-square distribution, IEEE Trans. Inform. Theory., 45(4), 1327-1333 (1999). DOI: 10.1109/18.761294

[8]. Annamalai, A., Tellambura, C., Matyjas, J.: A new twist on the generalized Marcum Q-function QM(a, b) with fractional-order M and its applications, Consumer Communications and Networking Conference, January 10-13: Las Vegas, NV, USA (2009).

[9]. Gil, A., Segura, J., Temme, N.M.: Computation of the Marcum Q-function, ArXiv, abs/1311.0681 (2013).

[10]. Olver, F.W.J., Lozier. D.W., Boisvert. R.F., Clark. C.W.: NIST Handbook of mathematical functions, Cambridge University Press (2010).

[11]. Bell, W.W.: Special Functions for Scientists and Engineers, Canada, D. Van Nostrand Company Ltd (1968).

[12]. Szegö, G.: Orthogonal Polynomials, Rhode Island, American mathematical society providence (1967).

[13]. Chatterjea, S.K.: On a generating function of Laguerre polynomials, Bollettino dell'Unione Matematica Italiana, Series 3, 17, 179-182 (1962).

[14]. Pent, M.: Orthogonal polynomial approach for the Marcum Q-function numerical computation, Electronic Lett. 4, 563-564 (1968). DOI: 10.1049/el:19680436

[15]. Andras, S., Baricz, A., Sun, Y.: The generalized Marcum Q-function: an orthogonal polynomial approach, Acta Universitatis Sapientiae, Mathematica, 3(1), 60-76 (2011). arXiv:1010.3348 [math.CA]